\documentclass[12pt,a4paper,reqno]{amsart}
\usepackage{amsmath,amsthm}
\usepackage{amsfonts}
\usepackage{amssymb}
\allowdisplaybreaks

\newcommand{\be}{\begin{equation}}
\newcommand{\ef}{\end{equation}}
\chardef\bslash=`\\ 

\newtheorem*{thm*}{Theorem}

\theoremstyle{definition}

\newtheorem*{remark*}{Remarks}
\newtheorem*{defn*}{Definition}

\theoremstyle{remark}

\newcommand{\G}{\Gamma}
\newcommand{\wt}{\widetilde}
\newcommand{\wh}{\widehat}
\newcommand{\fc}{\frac}

\newcommand{\iy}{\infty}
 \renewcommand{\sectionmark}[1]{}

\newcommand{\ve}{\varepsilon}

 \usepackage{amsfonts}
\newcommand{\field}[1]{\mathbb{#1}}
\newcommand{\g}{\gamma}

\newcommand{\dl}{\delta}
\newcommand{\D}{\field{D}}

\newcommand{\z}{\zeta}
\newcommand{\ov}{\overline}
\newcommand{\vp}{\varphi}
\newcommand{\hC}{\wh{\field{C}}}
\newcommand{\C}{\field{C}}
\newcommand{\R}{\field{R}}

\newcommand{\B}{\mathbf{B}}
\newcommand{\T}{\mathbf{T}}
\newcommand{\Belt}{\operatorname{Belt}}

\newcommand{\Fib}{\operatorname{Fib}}
\newcommand{\Teich}{\operatorname{Teich}}

\newcommand{\x} {\mathbf x}
\renewcommand{\a} {\alpha}

\begin{document}

\title{On evaluation of polynomial coefficients}
\author{Samuel L. Krushkal}

\begin{abstract} The evaluation of polynomial coefficients is an old mathematical problem, which is
still actively investigated by many authors.
The paper strengthens the known results in this direction and provides a complete solution for the case
when the space $P_n$ of polynomials is equipped with Hilbert's norm.

Our approach involves the intrinsic features of univalent functions and of Teichm\"{u}ller spaces.
\end{abstract}

\date{\today\hskip4mm({Eval.tex})}

\maketitle

\bigskip

{\small {\textbf {2020 Mathematics Subject Classification:} Primary: 30C10, 30F60;
Secondary 12D05,  30C50, 30C75, 31A05}

}

\bigskip
\markboth{Samuel L. Krushkal}{On evaluation of polynomial coefficients} \pagestyle{headings}

\bigskip\bigskip
\centerline{\bf 1. INTRODUCTORY REMARKS}.

\bigskip\noindent
{\bf 1.1}.
Estimating coefficients of polynomials and of holomorphic functions is one of
the oldest mathematical topics, but it still is actively investigated and remains an important problem in many algebraic, geometric and physical applications of complex analysis.

The following results were established recently in \cite{Kr6}, \cite{Kr7}.
Consider the collection $P_n$ of polynomials
$$
p_n(z) = c_0 + c_1 z + \dots + c_n z^n
$$
of degree $n \ge 1$, which is a linear space isomorphic to $C^{n+1}$, and regard $P_n$ as a subspace
of some Banach space $X$ of holomorphic functions on the unit disk $\D = \{|z| < 1\}$.

We regard the polynomials $p_n$ also as the {\bf Schwarzian derivatives}
$$
S_w(z) = \left(\frac{w^{\prime\prime}(z)}{w^\prime(z)}\right)^\prime
- \frac{1}{2} \left(\frac{w^{\prime\prime}(z)}{w^\prime(z)}\right)^2
$$
of locally univalent holomorphic functions $w(z) = a_0 + a_1 z+ \dots$ on $\D$ arising as the solutions of the differential equation
 \be\label{1}
(w^{\prime\prime}(z)/w^\prime(z))^\prime - (w^{\prime\prime}(z)/w^\prime(z))^2/2 = p_n(z)
\end{equation}
satisfying the prescribed initial conditions
$w(0) = a_0, \ \ w^\prime(0) = a_1, \ \ w^{\prime\prime}(0) = a_2$ (equivalently, as ratios of two linearly independent solutions $\eta_1, \eta_2$ of the linear differential equation
$2 \eta^{\prime\prime} + f(z) \ \eta = 0$
satisfying $\eta_1(0) = 0, \ \eta_1^\prime(0) = 1$ and $\eta_2(0) = 1, \ \eta_2^\prime(0) = 0$).

Generically, the solutions of such equations are connected with the special functions, cf. \cite{Ka}.

The Schwarzian derivatives are hyperbolically bounded holomorphic functions in $\D$, that means have the growth
$$
f(z) = O(1/(1 - |z|)^2) \quad \text{as} \ \ |z| \to 1.
$$
The complex Banach space $\B$ of all such functions with norm
$$
\|f\|_\B = \sup_\D (1 - |z|^2)^2 |f(z)|
$$
is dual to the space $A_1(\D)$ of integrable holomorphic functions on $\D$ (and $A_1(\D) \subset\B$).

Pick a rotationally invariant bounded convex subdomains $G \subset P_n$ (that is $G$  contains with every their polynomial $p_n$ its pre and post rotations
 \be\label{2}
p_n(z) \mapsto p_{n,\a, \beta}(z) =  e^{i \beta} p_n(e^{i \a} z)
\end{equation}
with independent $\a, \beta \in [- \pi, \pi]$. The Minkowski functional of any such domain $G$ generates a norm on $\C^{n+1}$ in which $G$ is the unit ball. Note also that all Banach norms on $P_n$  are equivalent.

One can distinguish in any bounded rotationally invariant convex domain $G \subset P_n$ two extremal polynomials determined up to rotations (2):
 \be\label{3}
p_n^0(z) = c_0^0 + c_1^0 z + \dots + c_n^0 z^n,
\end{equation}
which has the maximal first coefficient on $G$, i.e., with $|c_1^0| = \max_G |c_1(p_n)|$, and the Schwarzian
 \be\label{4}
S_{w_0}(z) = c_{0,0} + c_{m_1,0} z^{m_1} + \dots + c_{m_j,0} z^{m_j} + \dots + c_{n,0} z^n
\end{equation}
of univalent solution of the equation (1) with the maximal second coefficient
$$
a_2(w_0) = a_2(G) := \sup \big\{|a_2(w)|: \ S_w \in G \big\}.
$$
There are domains, for which these extremal polynomials are different. Such situation occurs
when $c_{1,0} = 0$.

\bigskip\noindent
{\bf 1.2}.
An interesting set in Banach spaces $X$ of holomorphic functions on the disk $\D$ is their subset $\mathcal Z_X$ formed by functions $f$ which are zero free on $\D$. Such functions arise and play an essential role in many fields of complex analysis, in particular, as the derivatives of locally univalent functions and of covering maps, and have been investigated by many authors.

The problem of estimating the coefficients of nonvanishing holomorphic functions in Banach spaces was stated by Hummel, Scheinberg and Zalcman in \cite{HSZ}. It is a special case of the evaluation problem
for coefficients of holomorphic functions on $\D$.

\bigskip\noindent
{\bf 1.3}.
The results of \cite{Kr6}, \cite{Kr7} mentioned above are presented by the following

\bigskip\noindent
{\bf Theorem 1}. {\it For any $n \ge 5$ and any rotationally invariant domain $G \subset P_n$ equipped with a norm $\|\cdot\|_X$, we have:

\bigskip\noindent
(a) \ $c_{1,0} \ne 0$.

\bigskip\noindent
(b) The extremal polynomials (3) and (4) for any such $n$ are equal (up to rotations (2));

\bigskip\noindent
(c) The coefficients $c_2, \dots, c_{n-1}$ of all polynomials $p_n \in G$  are sharply estimated by
$$
|c_m| \le \max  \ (|c_{1,0}|, |c_{m,0}|) = \max (|c_1^0|, |c_m^0|), \quad m = 2, \dots, n - 1.
$$

\bigskip\noindent
(d) The similar assertions and estimates are valid also for nonvanishing polynomials on $\D$, i.e., with zeros located outside of $\D$.   }

\bigskip
This theorem is valid for a wide collection of Banach norms on $P_n$, especially, when $P_n$ is embedded into the spaces $L_q, \ 2 \le q \le \iy$ or into a Hilbert space.
We need (in the proof) that these norms are compatible with the $\B$-norm.

As for the part $(d)$, an interesting set in any Banach space $X$ formed by holomorphic functions  is
its subset $\mathcal Z_X$ formed by the zero-free functions $f$. Such functions arise and play an essential role in many fields of complex analysis, in particular, as the derivatives of locally univalent functions and of covering maps.

\bigskip
The problem of estimating the coefficients of nonvanishing holomorphic functions $f$ from the general Banach spaces $X$ was stated by Hummel, Scheinberg and Zalcman in \cite{HSZ} in the form: find the value of
$$
c_n(\mathcal Z_X) = \sup_{\mathcal Z_X} \ \{|c_n(f)| : \ \|f\|_X \le 1\}.
$$
It is a special case of the evaluation problem for coefficients of holomorphic functions on $\D$.

The approach of \cite{Kr6}, \cite{Kr7} is more general; it concerns estimating on arbitrary convex domains in $P_n$  with a prescribed norm.

The Minkowski functional $l_G$ of a bounded rotationally invariant convex domain $G$ in a locally convex linear space $X$ generates a norm on this space, equivalent to the initial norm on $X$, and in the new norm this domain is the unit ball. This implies for $P_n$ that after the isometrical embedding of $P_n$ into this space $X$, the evaluation problem is reduced to estimating the quantities
$$
c_m(P_n) = \sup_G \ \{|c_m(p_n)| : \ \|p_n\|_X \le 1\}.
$$

\bigskip\bigskip
\centerline{\bf 2. MAIN RESULTS: CANONICAL FORMS OF EXTREMAL POLYNOMIALS}
\centerline{\bf FOR SOME HILBERT NORMS}

\bigskip\noindent
{\bf 2.1.Main theorems}.
We consider the Hilbert spaces $H$ of holomorphic functions
$$
f(z) = c_0 + c_1 z + c_2 z^2 + \dots + c_n z^n + \dots
$$
containing all powers $z^n, \ 0, 1, \dots$ and thereby the collections $P_n$ and are focused on the unit ball
$$
B_{1,H}(P_n) = \{p_n: \ \|p_n\|_H < 1\}.
$$
The simplest case arises when $P_n$ is equipped with the canonical $l^2$ norm
 \be\label{5}
\|p_n\|_{l^2} = \Bigl( \sum\limits_0^n |c_n|^2\Bigr)^{1/2}
\end{equation}
pulled back from $\C^{n+1}$.

Our goal is to show that in the Hilbertian case Theorem 1 can be essentially strengthened by giving the extremal polynomials with maximal $|c_m|$ explicitly.
This is provided by the following two theorems, which are new even in the case of norm (5).

\bigskip\noindent
{\bf Theorem 2}. {\it Let $n > 2$. Then for any $m = 2, 3, \dots, n - 1$ the extremal polynomial $p_{n,m}^0(z)$ maximizing the coefficient $c_m$ on the ball $B_{1,H}(P_n)$ is of the form
 \be\label{6}
p_m^0(z) =  c_1^0 z^m
\end{equation}
with
$$
|c_1^0| = \max_{B_{1,H}(P_n)} \ |c_1| = |S_{w_0}^\prime(0)|.
$$
The value $|c_1^0|$ is uniquely determined by the maximal modulus of the second coefficients  $a_2$
of univalent solutions $w_0(z)$ to the equation (1) among polynomials $p_n$ with $\|p_n\| \le 1$, up
to rotations (2). }

\bigskip
The assertion of this theorem requires some explanations.
One has to preserve by variations of $p_n$ two quantities: the maximal value $|c_m^0|$ and the
admissible maximal value of $\|p_n^0\|_H$ generating the univalent solutions of (1) on $\D$.

It turns out that one needs for this only the coefficients $c_m, c_n$ and that other coefficients
$c_s$ with $s \ne m, n$ do not affect $\max |c_m^0|$.
This is caused by the fact that we work on the maximal ball $B_{\rho,H} = \{\|p_n\|_H \le \rho\}$ touching the boundary of the space $\T$ from inside and therefore one cannot generically one vary arbitrarily the maximal values $|c_m^0|$ and $\|p_n^0\|_H$  without disturbing the univalence.

\bigskip
The next theorem concerns the nonvanishing (zero free) polynomials on $\D$.

\bigskip\noindent
{\bf Theorem 3}. {\it The extremal polynomial maximizing $|c_m|$ on the set $\mathcal Z_{P_n}$ of nonvanishing polynomials in the unit disk with maximal $|c_m|$ is (up to rotations (2)) of the form
 \be\label{7}
p_m^0(z) = \fc{c}{2} (z^m +1)
\end{equation}
with $c > 0$ uniquely determined by the maximal second coefficient $a_2$ of univalent solutions
$w(z)$ to the equation (1). }

\bigskip
Similar to \cite{Kr6}, \cite{ Kr7}, the proofs of these theorems in the general case essentially involve some deep geometric and analytic features of Teichm\"{u}ller spaces and of univalent functions.

\bigskip\noindent
{\bf 2.2. Remarks}.

\bigskip\noindent
{\bf 1}.
There is an exceptional case in Theorem 2, when its proof can be simplified.
Namely, if $H$ is the space $l^2$ with norm (5), the coefficients $c_m$ coincide with the Fourier coefficients and the assertion of Theorem 2 also follows by applying Parseval's equality.

\bigskip\noindent
{\bf 2}.
The upper bound for coefficient $a_2$ is given by Proposition 2 presented below.
It depends on the norm of $H$. The same holds for the constants $c_1$ and $c$ in (6), (7).

\bigskip\bigskip
\centerline{\bf 3. UNDERLYING RESULTS}

\bigskip\noindent
{\bf 3.1. Quasiconformal deformations associated with polynomials}.
One of the main underlying ingredients in the proofs of Theorems 2 and 3 is the following result from
\cite{Kr4}, on quasiconformal deformations of holomorphic functions ensuring the existence of quasiconformal maps $h$, which are conformal outside of a given domain $D \subset \hC$, have the prescribed values $h^{(\a_j)}$ of a finite collection of their derivatives in the prescribed points of $\ov D$ and preserve the norm of the original holomorphic function or, more generally, disturbing this norm in the controlled bounds.
The origins of such results go back to quasiconformal deformations of the punctured Riemann surfaces of finite conformal type, that is with the finitely generated fundamental groups, considered in \cite{Kr1}.

In the case of polynomials, these deformations imply the biholomorphic maps of small neighborhoods of polynomials with a prescribed distortion of norms.

Fix a natural $n >  1$ and consider the vectors $\mathbf d = (d_0, d_1, \dots, d_n) \in \C^{n+1}$ with the Euclidean norm $|\mathbf d|$ on $\C^{n+1}$.

\bigskip\noindent
{\bf Proposition 1}. \cite{Kr5} {\it Let $P_n$ be equipped with some Hilbert norm
$\| \cdot\|$, i.e., embedded into a Hilbert space $H$.  Then for any polynomial
$p_n(z) = c_0 + c_1 z + \dots + c_n z^n$ with $c_n \ne 0$, there exists a
number $\ve_0 > 0$ such that for every point
$$
{\mathbf d} = (d_0, d_1, \dots, \ , d_{n-1}) \in \C^n
$$
and every $a \in \R$ with
$$
|\mathbf  d| \le \ve, \ \ |a| \le \ve, \ \ \ve < \ve_0,
$$
there exists a quasiconformal automorphism $h$ of the extended complex plane $\hC, \ h(0) = 0$, which is conformal outside of some disk $\D(w_0, r) \subset p_n(\D^*)$ and determines a polynomial
$p_n^*(z) = \sum\limits_0^n \ c_m^* z^m$ with norm $\|p_n^*\|_H =  \|p_n \|_H + a$ and coefficients
$$
c_m^* = h^{(m)}(p_m(0)) = m! (c_m + d_m), \quad  m = 0, 1, \dots, n - 1.
$$
The Beltrami coefficient $\mu_h$ of $h$ satisfies $\|\mu_h \|_\iy \le M \ve$.
The quantities $\ve_0$ and $M$ depend only on $f$ and $n$}.

\bigskip
This proposition implies that the points $(h^\mu(c_1), \dots, h^\mu(c_{n-1}))$ fill a whole neighborhood
of the original point $(c_1, \dots, c_{n-1})$ in $\C^{n-1}$, while the deformation of the last coordinate $c_n$ is used for preserving the prescribed norm on $P_n$.

Such result holds also for the weighted spaces $L_q$, Hardy spaces $H^q$ with $2 \le q < \iy$, etc.;
see \cite{Kr2}.

Let us also mention, in connection with Theorem 3,  that the variations $h^\mu$ created by
Proposition 1 preserve nonvanishing of functions, i.e., $f$ and $h^\mu \circ f$ are zero free
simultaneously.

\bigskip\noindent
{\bf 3.2. Covering estimates}. We also will essentially apply the following covering result, which
generalizes Koebe's one-quarter theorem.

Let $G$ be an open connected rotationally invariant submanifold in the universal Teichm\"{u}ller space $\T$, which is a biholomorphic image of a domain in Banach space $X$.

\bigskip\noindent
{\bf Proposition 2}. \cite{Kr3} {\it Let $w(z)$ be a univalent solution of the equation (1) satisfying $w(0) = 0, \ w^\prime(0) = e^{i \theta}$ with the fixed $\theta \in [-\pi, \pi]$.
and let $w_0(z) = e^{i \theta} z + a_2^0 z^2 + \dots$ be one of the maximizing functions for $a_2$. Then the image domain $w(\D)$ covers entirely the disk
$D_{1/(2 |a_2(G)|)} = \{|w| < 1/(2 |a_2(G)|)\}$.
The radius value $1/(2 |a_2(G)|)$ is sharp for this collection of functions, and the circle $\{|w| = 1/(2 |a_2(G)|)$ contains the points not belonging to $w(\D)$ if and only if $|a_2| = |a_2(G)|$
(i.e., when $w$ is one of the maximizing functions for $a_2$ on the set $G$).

The inverted functions
$$
W(z) = 1/w(1/z) = e^{-i \theta} z - a_2^0 + b_1 z^{-1} + b_2 z^{-2} + \dots
$$
with $\z \in \D^*$ map the disk $\D^*$ onto a domain whose boundary is entirely contained in the
disk $\{|W + a_2(G)| \le |a_2(G)|\}$.
}

\bigskip\noindent
{\bf 3.3. A general distortion result for coefficient functionals}.
Let $X$ be a complex submanifold in the space $\B$. Consider the coefficient functionals
$$
I(f) = I \bigl(c_{m_1}, \dots, c_{m_l} \bigr)
$$
depending on a distinguished finite set of coefficients of functions $f \in X$ and take the univalent
functions $w(z) = z + a_2 z^2 + dots$ in unit disk with the Schwarzians $S_w = f$. Then the coefficients
$c_l$ are the polynomials from the initial cofficients $a_2, \dots, a_s$ of functions $w(z)$,
and hence, any indicated functional $I(f)$ generates the corresponding polynomial functional
$J(w) = J(a_2, \dots, a_s)$ on the whole class $S$.

The following Proposition describes estimating the general rotationally invariant coefficient
functionals in the norms of some Banach spaces $X$ (also including the Hilbert spaces and spaces $L_q$ with $q \ge 2$).

\bigskip\noindent
{\bf Proposition 3}. {\it Any rotationally invariant polynomial functional
$$
J(w) = J(a_2, \dots, a_s)
$$
on $S$ whose zero set $Z_J = \{w \in S: \ J(w) = 0\}$ is separated from the rotation set of extremal function of Proposition 2, is maximized on the space $X$ only by the rotations of this function.

In other words, any extremal function of $J(w)$ must be simultaneously maximal for the second coefficient $a_2$ on this class, unless $J(w_0) = 0$. }

\bigskip
The proof of this Proposition is given in \cite{Kr4} and in \cite{Kr5}. It essentially involves Proposition 2.

\bigskip\noindent
{\bf 3.4. Univalence}.
We also will use the relation between the $l^2$ and $\B$ norms on $P_n$. It implies the univalence
of solutions $w(z)$ to the equation (1) for $p_n \in B_1(P_n)$ on the entire disk $\D$ and quasiconformal
extension.

Identifying the polynomials $p_n$ with the points
$$
\mathbf c = (c_0, c_1, \dots, c_n, 0, \dots, 0, \dots) \in l^2,
$$
consider the balls
$$
B_r(P_n) = \{p_n: \ \|p_n\|_{l^2} < r\}, \quad 0 < r < \iy.
$$
Noting that all $p_n$ also belong to the Bergman spaces  $A_p(\D), p > 2$,  of holomorphic functions in the disk $\D$ with norm
$$
\|f\|_{A_p} = \Bigl( \fc{1}{\pi} \ \iint_\D |f(z)| dx dy \Bigr)^{1/p} \quad (z = x + i y).
$$
one obtains, applying the H\"{o}lder inequality,
$$
\|p_n\|_{A_1} = \fc{1}{\pi} \ \iint_\D |p_n(z)| dx dy \le \|p_n\|_{A_2} =\|p_n\|_{l^2}.
$$
Combining this with the well-known relations $A(\D) \subset \B$
$$
\|f\|_\B \le \|f\|_{A(\D)} \quad \text{if} \ \ f \in A(\D)
$$
(see, e.g., \cite{Be1}), one concludes that all $p_n$ from the ball $B_{1,H}(P_n)$ are the points of the
universal Teichm\"{u}ller space $\T$, and hence {\it such polynomials (more precisely, the corresponding quadratic differentials $p_n(z) dz^2$) are the Schwarzian derivatives of univalent functions (solutions to the equation (1)) in} $\D$, and these solutions adtmit quasiconformal extensions across the unit circle onto the extended complex plane $\hC$.

In fact, the Ahlfors-Weill theorem  yields that this is valid for $f$ from the bigger balls.
We use this in the proofs of the above theorems, which involve the maximal ball $B_{r_0,H}(P_n)$ of radius $r_0 > 1$, which touch from inside the boundary of the universal Teichm\"{u}ller space $\T$ (the  maximality provides an obstruction to quasiconformal extendability of univalent solutions to (1) corresponding to the common boundary  points).

Note also that since the balls in any normed space are convex sets, Theorems 2 and 3 also can be reformulated in the spirit of Theorem 1 in terms of convex sets with respect to the canonical $l^2$- metric mentioned in Section {\bf 2}.

\bigskip\bigskip
\centerline{\bf 4. REMARKS ON TEICHM\"{U}LLER SPACES}

\bigskip
The proofs of Proposition 3 (underlying Theorem 2) and of Theorem 3 involve essentially some deep geometric and analytic features of Teichm\"{u}ller spaces. To make the exposition self-contained, we present the needed results on these spaces; for details see, e.g., \cite{Be}, \cite{GL}, \cite{Le}.

It is technically more convenient to deal with the inverted functions $W(z) = 1/w(1/z)$ which are univalent (and zero free) on the disk $\D^*$ and admit quasiconformal extension to onto $\ov \D$).
Such univalent functions in $\D^*$ form the class $\Sigma$.

\bigskip\noindent
The {\bf universal Teichm\"{u}ller space} $\T = \Teich(\D)$ is the space of quasisymmetric homeomorphisms of the unit circle $\mathbb S^1 = \partial \D$ factorized by M\"{o}bius maps;  all Teichm\"{u}ller spaces have their biholomorphic copies in $\T$.

The canonical complex Banach structure on $\T$ is defined by factorization of the ball of the Beltrami coefficients (or complex dilatations)
$$
\Belt(\D)_1 = \{\mu \in L_\iy(\C): \ \mu|\D^* = 0, \ \|\mu\| < 1\},
$$
letting $\mu_1, \mu_2 \in \Belt(\D)_1$ be equivalent if the corresponding  quasiconformal maps $W^{\mu_1}, W^{\mu_2}$ (solutions to the Beltrami equation $\partial_{\ov{z}} W = \mu \partial_z W$
with $\mu = \mu_1, \mu_2$) coincide on the unit circle $\mathbb S^1$ (hence, on $\ov{\D^*}$). Such $\mu$ and the corresponding maps $w^\mu$ are called $\T$-{\it equivalent}.

We shall use the following lemma, which allows one to use some other normalizations of
quasiconformally extendable functions.

\bigskip\noindent
{\bf Lemma 1}. \cite{Kr4} {\it For any Beltrami coefficient $\mu \in \Belt(\D^*)_1$ and any $\theta_0 \in [0, 2 \pi]$, there exists a point $z_0 = e^{i \a}$ located on $\mathbb S^1$ so that
$|e^{i \theta_0} - e^{i \a}| < 1$ and such that for any $\theta$ satisfying
$|e^{i \theta} - e^{i \a}| < 1$ the equation
$\partial_{\ov z} w =  \mu(z) \partial_z w$
has a unique homeomorphic solution $w = w^\mu(z)$, which is holomorphic on the unit disk $\D$
and satisfies
$$
w(0) = 0, \quad w^\prime(0) = e^{i \theta}, \quad w(z_0) = z_0.
$$
Hence, $w^\mu(z)$ is conformal and does not have a pole in $\D$ \ (so
$w^\mu(z_{*}) = \iy$ at some point $z_{*}$ with $|z_{*}| \ge 1$).  }

\bigskip
In particular, this lemma allows one to define the Teichm\"{u}ller spaces using the quasiconformally extendible  univalent functions $w(z)$ in the unit disk $\D$ normalizing these functions by
$$
w(0) = 0, \quad w^\prime(0) = e^{i \theta}, \quad w(1) = 1.
$$
All such functions are holomorphic in the disk $\D$.

\bigskip
The Teichm\"{u}ller space $\T_1 = \Teich(\D_{*})$ {\bf of the punctured disk} $\D_{*} = \D \setminus \{0\}$ is formed by classes $[\mu]_{\T_1}$ of $\T_1$-{\bf equivalent} Beltrami coefficients $\mu \in \Belt(\D)_1$ so that the corresponding quasiconformal automorphisms $w^\mu$ of the unit disk coincide on both boundary components (unit circle $\mathbb S^1$ and the puncture $z = 0$) and are homotopic on $\D \setminus \{0\}$.
This space can be endowed with a canonical complex structure of a complex Banach manifold
and embedded into $\T$ using uniformization of $\D_{*}$ by a cyclic parabolic Fuchsian
group acting discontinuously on $\D$ and $\D^*$. The functions $\mu \in L_\iy(\D)$ are lifted to
$\D$ as the Beltrami measurable $(-1, 1)$-forms  $\wt \mu d\ov{z}/dz$ in $\D$ with respect to
$\G$, i.e., via $(\wt \mu \circ \g) \ov{\g^\prime}/\g^\prime = \wt \mu, \ \g \in \G$,
forming the Banach space $L_\iy(\D, \G)$; we extend these $\wt \mu$ by zero to $\D^*$. Then $\T_1$ is canonically isomorphic to the subspace $\T(\G) = \T \cap \B(\G)$,
where $\B(\G)$ consists of elements $\vp \in \B$ satisfying $(\vp \circ \g) (\g^\prime)^2 = \vp$
in $\D^*$ for all $\g \in \G$.

Due to {\bf the Bers isomorphism theorem}, {\it the space $\T_1$ is biholomorphically isomorphic
to the Bers fiber space
$$
\Fib(\T) = \{(\phi_\T(\mu), z) \in \T \times \C: \ \mu \in \Belt(\D)_1, \ z \in W^\mu(\D)\}
$$
over the universal Teichm\"{u}ller space $\T$ with holomorphic projection $\pi(\psi, z) = \psi$} (see \cite{Be2}).

This fiber space is a bounded hyperbolic domain in $\B_2 \times \C$ and represents the collection of domains $D_\mu = w^\mu(\D)$ as a holomorphic family over the space $\T$. For every $z \in \D$,  its
orbit $w^\mu(z)$ in $\T_1$ is a holomorphic curve over $\T$.

The indicated isomorphism between $\T_1$ and $\Fib(\T)$ is induced by the inclusion map \linebreak
$j: \ \D_{*} \hookrightarrow \D$ forgetting the puncture at the origin via
$$
\mu \mapsto (S_{w^{\mu_1}}, w^{\mu_1}(0)) \quad \text{with} \ \
\mu_1 = j_{*} \mu := (\mu \circ j_0) \ov{j_0^\prime}/j_0^\prime,
$$
where $j_0$ is the lift of $j$ to $\D$.

The Bers theorem is valid for Teichm\"{u}ller spaces $\T(X_0 \setminus \{x_0\})$ of all punctured hyperbolic Riemann surfaces $X_0 \setminus \{x_0\}$; we use only its special case.

\bigskip
The spaces $\T$ and $\T_1$ can be weakly (in the topology generated by the spherical metric on $\hC$) approximate by finite dimensional Teichm\"{u}ller spaces $\T(0, n)$ of punctured spheres (Riemann surfaces of genus zero)
$$
X_{\mathbf z} = \hC \setminus \{0, 1, z_1 \dots, z_{n-3}, \iy\}
$$
defined by ordered $n$-tuples $\mathbf z = (0, 1, z_1, \dots, z_{n-3}, \iy), \ n > 4$ with distinct
$z_j \in \C \setminus \{0, 1\}$ (the details see, e.g., in \cite{Kr3}).

Another canonical model of $\T(0, n) $ is obtained again using the uniformization. This space
is biholomorphic to a bounded domain in the complex Euclidean space $\C^{n-3}$.

Note also that all Teichm\"{u}ller spaces are complete metric spaces with intrinsic Teichm\"{u}ller metric defined by quasiconformal maps. By the Royden-Gardiner theorem, this metric equals the hyperbolic  Kobayashi metric determined by the complex structure (see, e.g., \cite{EKK}, \cite{GL}, \cite{Ro}).

\bigskip\bigskip
\centerline{\bf 5. PROOF OF THEOREM 2}

\bigskip
First we show that for all $n > 2$ the extremal Schwarzian $S_{w_0}$ has the first coefficient
 \be\label{8}
c_{1,0} \ne 0,
\end{equation}
which implies by Proposition 3 the equality of both extremal polynomials $S_{w_0}$ and $p_n^0$ (up to rotations). This is one of the main points of the proof.

Using the linearity of $c_m(p_n)$ on $P_n$, we consider the maximal ball $B_\rho(P_n)$ of radius
$\rho \ge 2$ centered at the origin chosen so that $B_\rho(P_n) \subset \T$ and its boundary sphere
has common points with the boundary of $\T$. Note here that  by the Ahlfors-Weill theorem \cite{AW}, any function (quadratic differential) $f \in \B$ with norm $\|f\|_\B  = k < 2$ is the Schwarzian derivative $S_w = f$ of a function $w$ which is univalent on the disk $\D$ and admits $k$-quasiconformal extension across the unit circle $\mathbb S^1$ to $\hC$ with Beltrami coefficient canonically connected with $f$.

\bigskip
We establish nonvanishing of $c_{1,0}$ by applying the Grunsky coefficients, which are defined for functions
$W(z) = z + b_0 + b_1 z^{-1} + \dots \in \Sigma$ by the expansion
 \be\label{9}
\log \fc{W(z) - W(\z)}{z - \z} = - \sum\limits_{m, n = 1}^\iy \a_{m n} z^{-m} \z^{-n},
\quad (z, \z) \in (\D^*)^2,
\end{equation}
and satisfy the inequality
 \be\label{10}
\Big\vert \sum\limits_{m,n=1}^\iy \sqrt{m n} \ \a_{m n} \ x_m x_n \Big\vert \le 1
\end{equation}
for any sequence $\mathbf x = (x_n)$ from  the unit sphere $S(l^2)$ of the Hilbert space $l^2$ with
norm $\|\x\| = \bigl( \sum\limits_1^\iy |x_n|^2\bigr)^{1/2}$; hence $|\a_{m n}| \le 1/\sqrt{m n}$
for all $m, n$.
Here the principal branch of the logarithmic function is chosen. By the Grunsky theorem \cite{Gr}
(extended to multiply connected domains in \cite{Gr}, \cite{Mi}) the
inequality (10) is the necessary and sufficient condition for univalence of a holomorphic function on
the unit disk.

The Grunsky coefficients of functions $w(z) \in S$ are defined in the similar fashion; moreover,
the functions $w(z) \in S$ and $W(z) = 1/w(1/z) \in \Sigma$ have the same Grunsky coefficients. Letting
$$
g(z, \zeta) = \log[(w(z) - w(\z))/(z - \z)]
$$
one obtains from  (10) that near any finite $\z$,
$$
\fc{\partial^2 g(z, \z)}{\partial z \partial \z} = \fc{1}{6} S_w(z) + O((z - \z)^2),
$$
and hence the Schwarzian of $f$ is given by
$$
S_w(z) = 6 \lim\limits_{\z\to z} \fc{\partial^2 g(z, \z)}{\partial z \partial \z}
= 6 \ \sum\limits_{m+n=2}^\iy \ \a_{m n} z^{m+n-2} = 6 \ \a_{1,1} + 24 \ \a_{2,1} z + \dots.
$$
A simple calculation yields
$$
\a_{2,1} = \a_{1,2} = a_4 - 3 a_2 a_3 + 2 a_2^3,
$$
which implies for polynomials $p_n \in P_n$ the equality
$$
b_1 = 24(a_4 - 3 a_2 a_3 + 2 a_2^3).
$$
Hence, the equality $b_1 = 0$ defines in $P_n$ the algebraic set whose preimage in $S$ has the equation
 \be\label{11}
a_4 - 3 a_2 a_3 +2  a_2^3 = 0
\end{equation}
(which embraces, in particular, all $w_0(z)$ corresponding to the even polynomials $p_n(z)$).

Now we assume that $n > 2$ and start with the case $m = 1$.

Take the solution $w_0(z)$ of the equation (1) defined by the extremal polynomial (4), i.e., with the maximal $|a_2|$ on $G$ and assume, in the contrary to (9), that $c_{1,0} = 0$, and hence,
$\max_G |c_1|$ is attained on some other polynomial
$$
p_n^* = S_{w_{*}} \ne p_n^0.
$$
Then, applying Proposition 1, one can construct a quasiconformal deformation $h^\mu(\z)$ with Beltrami coefficient $\mu$ supported in some disk $\D_{\z_0}(\rho) = \{|z - \z_0| < \rho\}$ sufficiently far from the origin so that first $n + 1$ the coefficients of $S_{h^\mu \circ w_{*}}$ satisfy
$$
\begin{aligned}
|c_1(S_{h^\mu \circ w_{*}}| &=
|a_4(h^\mu \circ p_n) - 2 a_2(h^\mu \circ p_n) a_3(h^\mu \circ p_n) + 2 a_2(h^\mu \circ p_n)^3| = O(\|\mu\|_\iy) > 0, \\
c_m(h^\mu \circ p_n) &= c_m(p_n) ]\quad \text{for} \ \ m = 0, 2, \dots, n -1,  \\
|c_n(h^\mu \circ p_n)| &= |c_n(p_n) - |O(\|\mu\|_\iy) < |c_n(p_n)|,
\end{aligned}
$$
and the partial sum $p_{n,h}^*$ of the first $n + 1$ terms of expansion  of the varied function $S_{h^\mu \circ w_{*}}$ has the same norm as $p_n^*$.
This contradicts the extremality of $|a_2^0|$ and of $|c_1^0|$ and proves (8). In addition, we have obtained that necessarily
$$
c_0^0 = 0, \quad c_2^0 = 0, \dots, \quad c_n^0 = 0.
$$
Otherwise, if some coefficient $c_j^0$ of $p_m^0$ with $j < m$ is nonzero, one can construct using Proposition 1 an appropriate variation $h^\mu$ such that it increases the modulus of $c_1^0$ and
decreases the modulus of $c_m^0$ so that the $H$-norm of the vector $(c_0^0, c_1^0, \dots, c_m^0)$
after this variation will be the same; hence, this vector (point) will be located in the original
closed ball $\ov{B_{1,H}(P_n)}$).

This contradicts to extremality of polynomial $p_n^0$. Hence, $p_n^0$ must be of the form (6) (with $n = 1$)\emph{} and satisfies the assumptions of the rotational invariance in Proposition 3. This proposition yields that, the first coefficient $c_1$ is maximized by the polynomial
$$
S_{w_0}(z) = p_1^0(z) = c_1^0 z,
$$
up to pre and post rotations (2).
The same arguments are valid for any $m = 0$

If $m = 2, \dots , n$, one must apply Proposition 3 to plurisubharmonic functionals
$$
I_m(p_n) = \max (|p_m(z)|, |\wt p_m(z)|),
$$
taking the truncated polynomial $\wt p_m(z) = c_m z^m$.

Since the map $p_m \mapsto \wt p_m$ is holomorphic on $P_n$, one can apply to $I_m$ the above arguments, which implies the sharp equalities
$$
\sup_{B_{1,H}(P_n)} \ I_m(p_n) \le \max (|c_1^0|, |c_m^0|).
$$
Since the additional coefficients $c_s$ with $s \ne m$ must be zero, the maximal value of all $|c_m^0|$ coincides with $|c_1^0| = \max_{B_{1,H}(P_n)} \ |c_1|$, and the assertion of Theorem 2 follows.

\bigskip
Note that for $m > 1$,
$$
S_{w_0}(z) = p_m^0(z) = c_m^0 z^m.
$$

\bigskip\bigskip
\centerline{\bf 6. PROOF OF THEOREM 3}

\bigskip
We must check that all steps in the proof of Theorem 2 (including the underlying propositions) are also valid for nonvanishing polynomials.

First, we pass to polynomials $p_{n,r}(z) = \fc{1}{r} \ p_n(r z)$ with $r < 1$ close to  $1$ and consider their collection
$$
B_{1,H,r}(P_n) = \{p_{n,r}(z): \ p_n \in B_{1,H}(P_n)\}.
$$
Every $p_{n,r}$ has a neighborhood $U(p_{n,r}, \epsilon(r))$ in $B_{1,H,r}(P_n)$ filled by polynomials which are zero free in the disk $\D$.
The union
$$
B_{1,\mathcal Z,r}(P_n) = \bigcup_{p_n \in \mathcal Z} U(p_{n,r}, \epsilon(r))
$$
of such neighborhoods is a rotationally invariant domain.
Then for each $r$ we have
$$
\max_{\ov \D} |p_n (r z)| \ge \dl(r) > 0.
$$
This insures that a sufficiently small neighborhood of a nonvanishing $p_n^0$ is filled by the coefficients of polynomials also nonvaninshing in $\D$.

It suffices to deal with this domain using instead of $a_{2,n}$ the maximal value of $|a_2(w)|$ on the image of the set $B_{1,\mathcal Z,r}(P_n)$ in $S$, which provides the corresponding subharmonic functionals $|c_m(p_n)|$ on $B_{1,\mathcal Z,r}(P_n)$. Denote this maximal value by $|c_{2,n}^0|$.

\bigskip
Now we lift the coefficients onto the Teichm\"{u}ller spaces $\T$ and $\Fib(\T)$ as follows.
Consider the classes $S_{z_0, \theta}$ of univalent functions $w(z)$ on the unit disk $\D$ with expansions
$$
w(z) = a_1 z^{-1} + a_2 z^2 + \dots
$$
with
$a_1 = e^{i \theta}, \ - \pi < \theta \le \pi$, having the fix point at $z_0$ on the boundary unit circle  $\mathbb S^1$ and
admitting quasiconformal extensions across $\mathbb S^1$ to the whole Riemann sphere $\hC$ and take their union
$$
\wh S = \bigcup_{z_0 \in \mathbb S^1, \theta \in (-\pi,\pi]} S_{z_0,\theta}.
$$
The Beltrami coefficients $\mu_f(z) = \partial_{\ov z} f/\partial_z f$ of these extensions run over the unit ball
$$
\Belt(D^*)_1 = \{\mu \in L_\iy(\C): \ \mu(z)|D = 0, \ \ \|\mu\|_\iy  < 1\}.
$$

Now pass to the inverted functions $W_w(z) = 1/w(1/z)$ which form the corresponding classes $\Sigma_{z_0,\theta}$ of nonvanishing univalent functions on the disk $\D^*$ with expansions
$$
W(z) =  e^{- i \theta} z + b_0 + b_1 z^{-1} + b_2 z^{-2} + \dots, \quad  F(z_0) = z_0,
$$
and let
$\Sigma^0 = \bigcup_{z_0,\theta} \Sigma_{z_0,\theta}$.

The coefficients $a_n$ of $f(z)$ and the corresponding coefficients $b_j$ of $F_f(z)$ are related by
$$
b_0 + e^{2i \theta} a_2 = 0, \quad b_n + \sum \limits_{j=1}^{n}
\epsilon_{n,j}  b_{n-j} a_{j+1} + \epsilon_{n+2,0} a_{n+2} = 0,
\quad n = 1, 2, ... \ ,
$$
where $\epsilon_{n,j}$ are the entire powers of $e^{i \theta}$. This
successively implies the representations of $a_n$ by $b_j$ via
 \be\label{12}
a_n = (- 1)^{n-1} \epsilon_{n-1,0}  b_0^{n-1} - (- 1)^{n-1} (n - 2)
\epsilon_{1,n-3} b_1 b_0^{n-3} + \text{lower terms with respect to} \ b_0.
\end{equation}
This transforms any coefficient functional $J(w)$ on $S$ depending on a finite set of  distinguished
coefficients $a_{m_1}, \dots, a_{m_s}$  into a coefficient functional $\wt J(W^\mu)$
on $\Sigma^0$ depending on the corresponding coefficients $b_j$ and lifts the functionals  $\wt J(W)$ and $J(w)$ holomorphically onto the universal Teichm\"{u}ller space $\T$ modelled via bounded domain in the space $\B$ of Schwarzians $S_W^\mu$. Our functional $J(w) = c_n$ is expressed in terms of these  coefficients.

To lift $J$ onto the covering space $\T_1$, we again pass to functional $\wh J(\mu) = \wt J(W^\mu)$
lifting $J$ onto the ball $\Belt(\D)_1$ and apply the $\T_1$-equivalence, i.e., the quotient map
$$
\phi_{\T_1}: \ \Belt(\D)_1 \to \T_1, \quad \mu \to [\mu]_{\T_1}.
$$
Thereby the functional $\wt J(W^\mu)$ is pushed down to a bounded holomorphic functional
$\mathcal J$ on the space $\T_1$ with the same range domain.
Regarding the points of $\T_1$ as the pairs
$X_{W^\mu} = (S_{W^\mu}, F^\mu(0))$, where $\mu \in \Belt(\D)_1$ obey $\T_1$-equivalence, and using the relations (12), one comes to a holomorphic functional
$$
\mathcal J(X_{W^\mu}) = \mathcal J(S_{W^\mu}, \ t), \quad t = W^\mu(0).
$$
on $\T_1 = \Fib(\T)$ and have to investigate its restriction to the image in $\Fib (\T)$ of the original submanifold corresponding to $G$.

By Propositions 2, the boundary of any domain $W^\mu(\D^*)$ for any $W^\mu(z) \in \Sigma^0$ is located in the disk $\{|W - b_0| \le a_{2,n} \}$, and the second coordinate $t$ runs over some subdomain $D_\theta$ in the disk $\D_{2 a_{2,}} = \{|t| < 2 a_{2,n} \}$ containing the origin.

Finally, we have to establish the range domain of $W^\mu(0)$ for $S_{W^\mu}$ running over
$G_{\mathcal Z,r}$ and describe the boundary points of this domain.

The features of construction of the Bers fiber space $\Fib(\T)$ imply that image of $G_{\mathcal Z,r}$ in the space $\Fib(\T)$ also is a connected submanifold covering $G_{\mathcal Z,r}$.

We select a dense subsequence $\{\theta_1, \theta_2, \dots\} \subset [-\pi, \pi]$, getting the corresponding coefficient functionals on the classes $S_{z_0,\theta_j}$ and $\Sigma_{z_0,\theta_j}$, and consider the sequence of increasing products of the quotient spaces
 \be\label{13}
\mathcal T_s = \prod_{j=1}^s \ \wh \Sigma_{\theta_j}/\thicksim \
= \prod_{j=1}^s \{(S_{W_{\theta_j}}, W_{\theta_j}^{\mu_j}(0)) \} \ \simeq \T_1 \times \dots \times \T_1,
\end{equation}
where the equivalence relation $\thicksim$ again means $\T_1$-equivalence.
The Beltrami coefficients  $\mu_j \in \Belt(\D)_1$ are chosen here independently. For any $\T_1$, presented in the right-hand side of (13), the corresponding values of $W_{\theta_j}^{\mu_j}(0)$ run over some domain $D_{\a_j} \subset \C$, and the corresponding collection
$\beta = (\beta_1, \dots, \beta_s)$
of the Bers isomorphisms
$\beta_j: \ \{(S_{W_{\theta_j}}, W_{\theta_j}^{\mu_j}(0))\} \to \Fib(\T)$
determines a holomorphic surjection of the space $\mathcal T_m$ onto the product of $m$ spaces
$\Fib(\T)$. This provides the corresponding holomorphic maps
$$
\mathbf J = (\mathcal J(S_{W_{\theta_1}}), \dots, \mathcal J(S_{W_{\theta_s}}))
$$
whose (polydisk) norm satisfies
$$
\max_{B_{1,\mathcal Z,r}(P_n)^s} \|\mathbf J\| = \max_j |J(S_{W_{\theta_1}})| = \max_{B_{1,\mathcal Z,r}(P_n)} |c_m|.
$$

Now pick an extremal polynomial $p_m^0 \in B_{1,\mathcal Z,r}(P_n)$ with
$$
|\mathcal J(S_{W^\mu}, t)| = |c_m^0| = \max_{B_{1,\mathcal Z,r}(P_n)} |c_m(p_n)|.
$$
and consider the corresponding equation (1) for this $p_m^0$ and the corresponding univalent
solution $w_0(z)$.

Using the features of spaces $\T$ and $\T_1$ indicated in Section 4, one delivers that the upper semicontinuous regularization $u_r(t)$ of the function
$$
\sup_{S_{W^\mu}} \mathcal J(S_{W^\mu}, t)
$$
with the supremum over all $S_{W^\mu} \in \T$ admissible for a given $t = W^\mu(0) \in D_{\a_n}$
is subharmonic on the disk $\D_\a$ of some radius $\a \le 2 |a_{2,n}|$ centered at the origin;
hence, it attains the maximal value on the boundary circle (the details are given, for example, in \cite{Kr3}, \cite{Kr4}, \cite{Kr6}). This maximal function is circularly symmetric and positive on any circle $\{|t| = \rho\}$.

Since this construction involves the polynomials and univalent $w(z)$ with $a_2$ arbitrarily close to
$a_{2,n}$, we have its subharmonicity in the whole disk $\D_{2 |a_{2,n}|}$. So, for every $c_m$,
$$
\max_{\D_{2 |a_{2,n}|}} \ u_r(t) =  \sup_{B_{1,\mathcal Z,r}(P_n)} \ \mathcal |J(S_{W^\mu}, t)| =
\max_{B_{1,\mathcal Z,r}(P_n)} \ |c_m|,
$$
provided that $c_{m,0} \ne 0, \ m \ge 1$.

Since the corresponding extremal polynomials $p_{n,r}^0$ do not vanish in $\D$, their free
terms $c_{0,r} \ne 0$, while similar to the proof of Theorem 2, one obtains maximizing $c_{m,r}$ that
the coefficients $c_{s,r}^0$ different from $c_{0,r}$ and $c_{m,r}$ must vanish.

Finally, the assumption $p_{n,r}(z) \ne 0$ on $\D$ yields that the zeros of these polynomials have the moduli $|z_{j,r}| = |c_{0,r}|/|c_{1,r}| \ge 1$. Hence the maximal admissible value of $|c_{1,r}|$ is attained when
$$
|c_{1,r}| = |c_{0,r}| = r/2.
$$
This yields for the approximating polynomials $p_{n,r}$ the bounds and extremal polynomials of type (7) on the sets $B_{1,\mathcal Z,r}(P_n)$. These bounds are continuous in $r$ and uniformly estimated for all $r \le 1$.

Going to the limit as $r \to 1$, one obtains the desired form (7) of extremal nonvanishing polynomials on
$\mathcal Z_{P_n}$ completing  the proof of Theorem 3.

\bigskip\bigskip
\centerline{\bf 7. GENERALIZATION}

\bigskip
The main Theorems 2 and 3 can be extended to more general Banach spaces, for example, to the reflective Banach spaces. The details will be presented elsewhere.

\bigskip
{\small\it{ \leftline{Department of Mathematics, Bar-Ilan University, 5290002 Ramat-Gan, Israel} \leftline{and Department of Mathematics, University of Virginia,  Charlottesville, VA 22904-4137, USA}}

\end{document}